\documentclass[dvips,a4paper]{article}

\usepackage{amsmath}
\usepackage{amssymb}
\usepackage[english]{babel}
\usepackage[latin1]{inputenc}
\usepackage{color}
\usepackage{graphicx}
\usepackage{xspace}
\usepackage{enumerate}
\usepackage{hyperref}

\newtheorem{theor+}{Theorem}[section]
\newtheorem{propo+}[theor+]{Proposition}
\newtheorem{estim+}[theor+]{Estimate}
\newtheorem{imple+}[theor+]{Implementation remark}
\newtheorem{corol+}[theor+]{Corollary}
\newtheorem{notat+}[theor+]{Notation}
\newtheorem{lemma+}[theor+]{Lemma}
\newtheorem{obser+}[theor+]{Observation}
\newtheorem{defin+}[theor+]{Definition}
\newtheorem{examp+}[theor+]{Example}
\newtheorem{remar+}[theor+]{Remark}
\newtheorem{conje+}[theor+]{Conjecture}
\newtheorem{quest+}[theor+]{Question}
\newtheorem{pictu+}[theor+]{Picture}
\newtheorem{notes+}[theor+]{Notes}
\newtheorem{exerc+}[theor+]{Exercise}
\newtheorem{backg+}[theor+]{Background}
\newtheorem{fact+}[theor+]{Fact}
\newtheorem{abort+}[theor+]{Terminal state}

\newenvironment{theor}[1]{\begin{theor+}\label{#1}\slshape}{\end{theor+}}
\newenvironment{propo}[1]{\begin{propo+}\label{#1}\slshape}{\end{propo+}}

\newenvironment{lemma}[1]{\begin{lemma+}\label{#1}\slshape}{\end{lemma+}}

\newenvironment{defin}[1]{\begin{defin+}\label{#1}\upshape}{\end{defin+}}

\newenvironment{examp}[1]{\begin{examp+}\label{#1}\slshape}{\end{examp+}}

\newenvironment{abort}[1]{\begin{abort+}\label{#1}\slshape}{\end{abort+}}
\newenvironment{demo}{{\bf\sl Proof:}}{\hfill$\square$}

\newcommand{\mikkel}[1]{#1} %\mbox{\textsf{CENSUR!}}}

\newcommand{\eight}[2]{\overbrace{#2}^{#1}}
\newcommand{\sixteen}[2]{\overbrace{#2}^{#1}}
\newcommand{\eighto}[1]{\overbrace{0,\dots,0}^{#1}}
\newcommand{\sixteeno}[1]{\overbrace{0,\dots,0,0,\dots,0}^{#1}}
\newcommand{\lego}{\textsf{LEGO}\xspace}
\newcommand{\lbls}{\textsf{LEGO} blocks\xspace}
\newcommand{\lbl}{\textsf{LEGO} block\xspace}
\newcommand{\bls}{blocks\xspace}
\newcommand{\bl}{block\xspace}
\newcommand{\RR}{\mathbb{R}}
\newcommand{\NN}{\mathbb{N}}

\newcommand{\tbf}{$2\times4$\xspace}
\newcommand{\bbw}{\mbox{$b \times w$}\xspace}
\newcommand{\tbfm}{{2\times4}}
\newcommand{\bbwm}{{b\times w}}
\newcommand{\FAIL}{\mathsf{FAIL}}
\newcommand{\AAA}{\mathcal{A}}
\newcommand{\aaa}{\mathsf{a}}
\newcommand{\CCC}{\mathcal{C}}
\newcommand{\ccc}{\mathsf{c}}
\newcommand{\BBB}{\mathcal{B}}
\newcommand{\bbb}{\mathsf{b}}
\title{On the entropy of 
  \textsf{LEGO}$^{\mbox{\small\textregistered}}$\footnote{\lego is a
    trademark of \lego Company}}
\author{Bergfinnur Durhuus and Søren Eilers}
\date{March 2005}
\begin{document}
\maketitle
\section{Introduction}
It has long been asserted that the number of ways to combine six
\tbf \lbls of the same color is
\[
102981500
\]
This number was computed  at \lego in 1974 (\cite{jkc}) and has
been systematically repeated, for instance in
\cite[p. 15]{lcp} and \cite{ulb}. Consequently, the number can be found in several ``fun
fact'' books and on more than 250 pages on the World Wide Web\footnote{{\href{http://www.google.com/search?q=102981500&ie=UTF-8&oe=UTF-8}{Google search}},
  October 2004}.
However, this number only gives (with a small error, as we shall see) the number of ways to build a tower
of \lbls of height six. The total number of configurations is 
\[
\mikkel{915103765}
\]
as found by independent computer calculations by the second author and
by Abrahamsen \cite{ma:data}. This figure has now been accepted by
\lego Company, cf.\ \cite{llife}.

\begin{figure}[!ht]
\begin{center}
\includegraphics[width=12cm]{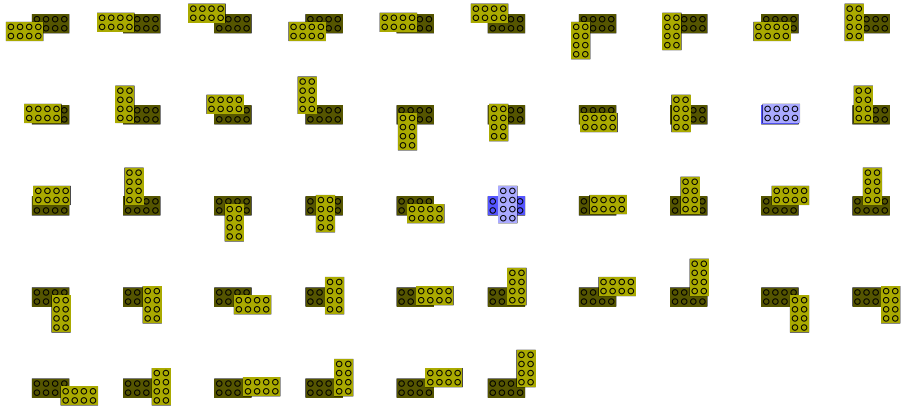}
\end{center}\caption{46 basic positions}\label{xxxxvi}
\end{figure}

We consider contiguous buildings of \lbls, disregarding color,  and identify them
up to translation and rotation. We think of a \bbw \lbl as a subset of
$\RR^3$ of
the form
\[
[x;x+b]\times [y,y+w]\times[z;z+1]
\]
or
\[
[x;x+w]\times [y,y+b]\times[z;z+1]
\]
Since the top and bottom of a  \lbl
are distinguishable we will only consider rotations in the
$xy$-plane.

It is easy to see that one such block may be put on top of
another in
\begin{equation}\label{oononsquare}
(2b-1)(2w-1)+(b+w-1)^2
\end{equation}
different ways if $b\not=w$ and in 
\begin{equation}\label{oosquare}
(2b-1)^2
\end{equation}
different ways if $b=w$.

With $w=4$ and $b=2$ we get $46$ possibilities, and note (as depicted
in blue on Figure \ref{xxxxvi}) that $2$ of these are
symmetric. Thus, letting $H_\tbfm(n,m)$ denote the number of ways to build a building of height $m$
with $n$ \tbf \lbls one then clearly has
\[
H_\tbfm(n,n)=\frac12(46^{n-1}+2^{n-1})
\]
Note that $H_\tbfm(6,6)=102981504$, so that in
fact \lego's computation is off by four.

By combining results of computer-aided enumerations with elementary
combinatorics one can further establish
\begin{equation}\label{oneless}
H_\tbfm(n,n-1)=46^{n-4}(-89115+37065n)+2^{n-4}(-8+16n)
\end{equation}
for $n\geq 3$ and 
\begin{gather}
H_\tbfm(n,n-2)=2^{n-7}(1785 - 825 n + 256 n)\label{twoless}\\
+46^{n-7}(-918674675 - 5330182078 n + 1373814225 n^2)\notag
\end{gather}
for $n\geq 5$, but as the problem is rather hopelessly non-markovian  there seems to
be no way to give formulae for the number of buildings of relatively low
height, or indeed for the total number $T_\tbfm(n)$ of  contiguous configurations,
counted up to symmetry. Although symmetry arguments and other tricks can be used to
prune the search trees somewhat, we are essentially left
with the very time-consuming option of going through all possible
configurations  to determine these numbers, which even with efficient computers seems completely
out of range for numbers such as $T_\tbfm(12)$. A sample of our results
may be seen in Figure \ref{obs}.

\begin{figure}[h]
\begin{center}
\begin{tabular}{|c||r|r|r|r|r|}\hline
$H_\tbfm(n,m)$&$m=2$&$m=3$&$m=4$&$m=5$&$m=6$\\\hline\hline
$n=2$&24&&&&\\\hline
$n=3$&500&1060&&&\\\hline
$n=4$&11707&59201&48672&&\\\hline
$n=5$&248688&3203175&4425804&2238736&\\\hline
$n=6$&7946227&162216127&359949655&282010252&102981504\\\hline
\end{tabular}
\end{center}
\caption{$H_\tbfm(n,m)$ for $m,n\leq 6$}\label{obs}
\end{figure}

\begin{figure}[h]
\begin{center}
\begin{parbox}{4.5cm}{\begin{tabular}{|c|r|}\hline
$n$&$T_\tbfm(n)$\\\hline\hline
1&1\\\hline
2&24\\\hline
3&1560\\\hline
4&119580\\\hline
5&10116403\\\hline
6&915103765\\\hline
7&85747377755\\\hline\end{tabular}}\end{parbox}
\qquad\begin{parbox}{6.2cm}{\includegraphics[width=6cm]{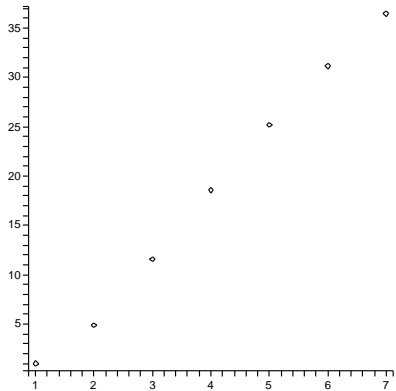}}
\end{parbox}\end{center}\caption{Rate of growth of $T_\tbfm$ with
semilogarithmic plot}\label{semiplot}
\end{figure}

These results seem to indicate, as shown on Figure \ref{semiplot}, that
$T_\tbfm(n)$ grows exponentially in $n$. In this paper we will show
that this is indeed the case, and give upper and lower bounds on the
rate of growth -- the \emph{entropy} of the blocks.

\section{Entropy of $\bbw$ blocks}

It is the goal of the present section to prove that the following
definition makes sense:

\begin{defin}{ent}
The \emph{entropy} of a \bbw \lbl is 
\begin{equation}\label{defeq}
s_{\bbwm}=
\lim_{n\longrightarrow\infty}\frac{\log(T_\bbwm(n))}{n}
\end{equation}
We let $h_{\bbwm}=\exp(s_\bbwm)$. 
\end{defin}

That the limit exists is by no means obvious, except of course for
$s_{1\times 1}=0$. We shall prove that this is the case in two steps,
first establishing convergence in $[0;\infty]$ and then proving that
the limit is finite.

It is inconvenient and irrelevant for our theoretical considerations
to identify buildings up to symmetry, so we establish definiteness in
another fashion. Suppressing the block size from the notation, we will
by $\AAA_n$ denote all contiguous
buildings containing $[0;w]\times [0;b]\times [0;1]$ with the further
property that there is no other block in $\RR^2\times [0;1]$ and no
block at all in $\RR^2\times [-1;0]$. Thus, the configuration can be
thought of as sitting on a base block at a fixed position.

We let $\aaa_n=\#\AAA_n$ and note

\begin{lemma}{sameentro}
We have
\[
\lim_{n\longrightarrow \infty}\frac{\log(T_\bbwm(n))}n=
\lim_{n\longrightarrow \infty}\frac{\log(\aaa_n)}n
\]
in the sense that if one limit exists, so does the other.
\end{lemma}
\begin{demo}
The claims follow immediately by the inequalities
\[
T_\bbwm(n-1)\leq\aaa_n\leq 4T_\bbwm(n)
\]
The leftmost inequality follows by mapping each equivalence class of
configurations with $n-1$ blocks to a representative placed on top of
$[0;w]\times [0;b]\times [0;1]$ and noting that this map is
injective. The rightmost follows by mapping each configuration to an
equivalence class and noting that this map is at most $4-1$.
\end{demo}

We now get

\begin{propo}{conv}
${\log(\aaa_n)}/n$ converges in $[0;\infty]$ as $n\longrightarrow
\infty$.
\end{propo}
\begin{demo}
One sees that $\aaa_{n+m}\geq \aaa_n\aaa_m$ by noting that an injective
map from $\AAA_n\times\AAA_m$ to $\AAA_{m+n}$ is defined by placing
the base block of the element of $\AAA_m$ somewhere on the top layer
of the element of $\AAA_n$.

Hence $\log(\aaa_n)$ is a superadditive sequence, and  $\log(\aaa_n)/n$
converges to $$\sup_{n\in\NN} \log(\aaa_n)/n.$$
\end{demo}

To prove that this limit is finite, i.e.\ that $\aaa_n$ grows no faster
than exponentially, we describe a surjective map
associating to each function
\begin{equation}\label{afunct}
S_n:\{1,\dots,2bw(n-2)\}\longrightarrow\{-bw,-bw+1,\dots,bw-1,bw\}
\end{equation}
an element of
\[
\AAA_n\cup\{\FAIL\}.
\]

Clearly the number of
such functions grows only exponentially in $n$. We shall subsequently look closer at which functions do indeed lead to
buildings with $n$ \lbls, and give much better estimates for $h_\bbwm$
than the obvious $(2bw+1)^{2bw}$.

With a fixed  enumeration of the studs and holes of a \bbw \lbl by the numbers 
$1,\dots,bw$, a map of the form \eqref{afunct}
gives rise to an element of $\AAA_n$, or the symbol $\FAIL$, as
follows.

 Take one
\lbl and call it block 1. Then read $S_n(1),\dots,S_n(bw)$ from left to
right to specify what to build
on top of block 1 as follows. If $S_n(1)>0$, take another \lbl and place
it parallely to block 1 with hole $S_n(1)$ on top of stud 1. If $S_n(1)<0$, 
take a \lbl and place
it orthogonally, rotated $+90^\circ$, to block 1 with hole $-S_n(1)$ on top of hole 1. In both cases,
give the new block the number 2. If $S_n(1)=0$,
do nothing. Then proceed to read $S_n(2)$ to see what, if anything, to
place on stud 2, and so on until $S_n(bw)$. Enumerate the blocks as they
are introduced.

\begin{abort}{abii}
If at any point a block collides with one which has
already been placed, the procedure terminates with $\FAIL$.
\end{abort}

These steps will result in the placing of between 0 and $bw$ blocks on
block 1.

\begin{abort}{abiii}
If at any point all $n$ blocks have been placed, consider the
unread values of $S_n$. If they are all $0$, the procedure
terminates successfully with an element of $\AAA_n$. If not, the
procedure terminates with $\FAIL$.
\end{abort}

\begin{abort}{abiv}
If, after reading the specifications for the first $m<n-1$ blocks, no
block $m+1$ has been introduced, the procedure terminates with $\FAIL$.
\end{abort}

We may now assume that a block 2 has been introduced and look at
$S_n(bw+1),\dots,S_n(2bw)$ which will specify what to build on top of  this block, if
anything, in the same way that $S_n(1),\dots,S_n(bw)$ specified what
to build on block 1. A positive number at $S_n(bw+1)$ will result in
the placing of a block on stud 1 of block 2
parallel to block 1, a negative number at $S_n(bw+1)$ will result in
the placing of a block on stud 1 of block 2 orthogonal to block
1, \emph{etc}. We proceed in the same way for blocks
$3,\dots,n-1$, but now read $2bw$ values where
$S_n((2m-4)bw+1),\dots,S_n((2m-3)bw)$ will specify what to put on top
of block $m$, and $S_n((2m-3)bw+1),\dots,S_n((2m-2)bw)$ will specify
what to put on underneath it 
in an analogous way.

\begin{abort}{abv}
If at any point a second  block is placed at the level $\RR^2\times [0;1]$, the procedure terminates with $\FAIL$.
\end{abort}

\begin{abort}{abvi}
If $S_n$ has been read to the end, consider the number of blocks
placed. If it is less than $n$, the procedure terminates with $\FAIL$. If not, 
it terminates successfully with an element of $\AAA_n$.
\end{abort}

We repeat this until one of the terminal states are reached.

If the procedure does not fail, it will result in a building of $n$ contiguous blocks, and clearly any such building may be constructed in
this way. Thus, the number of possibilites for maps $S_n$ dominates
$\aaa_n$, as desired, and we have:

\begin{theor}{first}
The limit in \eqref{defeq} exists for any block dimension $\bbw$.
\end{theor}

We can give general bounds of $h_\bbwm$, but as wee shall see below in
the case $w=4,b=2$, these bounds can in general be rather dramatically
improved.  

\begin{theor}{genest}
If $b\not=w$ we have
\begin{equation}
(2b-1)(2w-1)+(b+w-1)^2\leq h_\bbwm\leq \frac{(2bw)^{2bw+1}}{(2bw-1)^{2bw-1}}
\end{equation}
If $b=w$ we have
\begin{equation}
(2b-1)^2\leq h_\bbwm\leq\frac{(b^2)^{2b^2+1}}{(b^2-1)^{2b^2-1}}
\end{equation}
\end{theor}
\begin{demo}
By \eqref{oononsquare} we clearly have
\[
\aaa_n\geq((2b-1)(2w-1)+(b+w-1)^2)^{n-1}
\]
when $b\not=w$, and this gives the lower bound in that case. For the
upper bound, note that a function $S_n$ with nonzero entries at $m$
locations will yield $\FAIL$ if $m\not =n-1$. Hence we get
\[
\aaa_n\leq
\left(\begin{array}{c}2bw(n-2)\\n-1\end{array}\right)(2bw)^{n-1}
\]
By Stirling's formula we then get
\begin{eqnarray*}
\aaa_n&\leq&\frac{(2bw(n-2))!}{(n-1)!((2bw-1)n-4bw+1)!}(2bw)^{n-1}\\
&=&\frac{(2bwn)!(2bw)^{n}}{n!((2bw-1)n)!}O(1)\\
&=&\frac{(2bwn)^{2bwn}(2bw)^{n}}{n^n((2bw-1)n)^{(2bw-1)n}}O(n^{-1/2})\\
&=&\frac{(2bw)^{2bwn}(2bw)^{n}}{(2bw-1)^{(2bw-1)n}}O(n^{-1/2})\\
&=&\left(\frac{(2bw)^{2bw+1}}{(2bw-1)^{2bw-1}}\right)^nO(n^{-1/2})
\end{eqnarray*}
from which the claim follows directly.

The square case follows similarly by \eqref{oosquare} and by noting
that functions $S_n$ may be chosen non-negative. 
\end{demo}

\begin{examp}{varexx}
We enumerate the studs of a \tbf \lbl according to Figure \ref{numerate}.
\begin{figure}
\begin{center}\sf
1 2 3 4\\
\includegraphics[width=2.5cm]{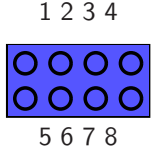}\\[-0.2cm]
5 6 7 8
\end{center}\caption{Enumeration of studs and holes}\label{numerate}
\end{figure}
Now consider functions $\{1,\dots,32\}\mapsto\{-8,\dots,8\}$ given by
\begin{gather}
(\eight{1}{0,5,0,0,-4,0,0,0},\eight{2}{0,0,0,0,-1,0,0,0},\sixteeno{3})\label{good}\\
(\eight{1}{0,-1,0,0,0,0,0,0},\eight{2}{0,5,0,0,0,0,0,0},\sixteen{3}{0,\dots,0,0,0,0,-1,0,0,0})\label{goodii}\\
(\eight{1}{1,1,0,0,0,0,0,0},\eighto{2},\sixteeno{3})\\
(\eight{1}{0,-1,0,0,0,0,0,0},\eight{2}{0,5,0,0,0,0,0,0},\sixteen{3}{0,\dots,0,0,0,0,-1,0,-1,0})\\
(\eighto{1},\eighto{2},\sixteeno{3})\\
(\eight{1}{0,-1,0,-1,0,0,0,0},\eighto{2},\sixteen{3}{0,0,2,0,0,0,0,0,\dots,0})\\
(\eight{1}{0,-1,0,0,0,0,0,0},\eight{2}{0,5,0,0,0,0,0,0},\sixteeno{3})
\end{gather}
where all ellipses indicate six consecutive zeros. The functions \eqref{good} and \eqref{goodii} give rise to the
buildings depicted on Figure \ref{goodbuildings}.
\begin{figure}
\begin{center}\sf
\includegraphics[width=2.5cm]{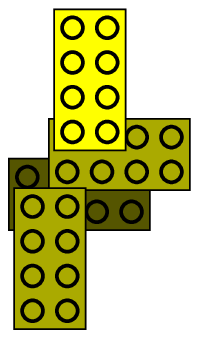}\qquad
\includegraphics[width=2.99cm]{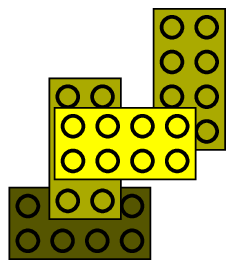}
\end{center}\caption{Buildings associated to \eqref{good} and \eqref{goodii}}
\label{goodbuildings}
\end{figure}
The remaining five functions give simple examples of functions
resulting in the procedure failing at Terminal state \ref{abii},
\ref{abiii}, \ref{abiv}, \ref{abv}, and \ref{abvi}, respectively.
\end{examp}

\section{Improved upper bounds}

In this section we shall describe methods to improve the upper bound on
$h_\bbwm$ given in Theorem \ref{genest}. They apply to any dimension,
but as they are somewhat \emph{ad hoc} we shall concentrate on our
favored dimension \tbf and leave other cases to the reader.

From Theorem \ref{genest} we know that $h_\tbfm\leq
16^{17}/15^{15}\leq 647.02$. We shall give a simple improved estimate leading 
to $h_\tbfm\leq 203.82$ and a somewhat more complicated one leading to 
$h_\tbfm\leq 191.35$. Besides being easier to state, the simpler estimate has
applications in producing statistical estimates for $\aaa_n$ for
relatively large
$n$.

Note that the surjective map associating buildings (or $\FAIL$) to certain maps
\[
S_n:\{1,\dots,16(n-2)\}\longrightarrow\{-8,-7,\dots,7,8\}
\]
is very far from being injective. We have already employed the fact
that unless the number of nonzero values is $n-1$, the function is
mapped to $\FAIL$. But we may also use that  the placement of a
block onto another may be indicated in $\ell$ different ways, where $\ell$ is
the number of studs of the lower block which are inserted into the upper block. Restricting attention to
maps where placements are indicated in a fixed way will not affect
surjectivity of the map.

Any partition of the 46 positions in Figure \ref{xxxxvi} into 8
sets $\mathcal{P}_1,\dots,\mathcal{P}_8$ with the property that any position in $\mathcal{P}_i$
employs stud $i$ of the lower block can be used to improve the upper
bound. One uses the convention that a position in $\mathcal{P}_i$ is always
indicated by a symbol at stud $i$, thus restricting 
the number of possibilities.

Another restriction is available when specifying what to add to block
$m$ for $m>2$. If we keep track of how block $m$ was introduced, we
know \emph{a priori} that one hole or one stud of it has already been
used, thus eliminating at least $16$ out of the $46$ possibilities on the
relevant side of the block. Dividing up the remaining 30 positions as
above, we get 64 sets $\mathcal{P}^j_i$ with the property that
$\mathcal{P}^i_i=\emptyset$ and that $\mathcal{P}^j_1,\dots,\mathcal{P}^j_8$ is a partition of the
30 positions which do not employ stud $j$.

\begin{theor}{sixest}
We have
\begin{eqnarray}\label{sixbounds}
\aaa_n\leq 
\left(\begin{array}{c}13n-23\\n-1\end{array}\right)6^{n-1}
\end{eqnarray}
and, consequently, that $h_\tbfm\leq 6\cdot13^{13}/12^{12}< 203.82$
\end{theor}
\begin{demo}
We partition the 46 positions into 8 sets, each consisting of 6 or 5
configurations, as indicated by the rows of Figure \ref{oneononeii}. Thus on
indices related to one side of blocks $1,\dots,n-1$, we need only allow
for 6 different symbols.

 \begin{figure}
\begin{center}
\includegraphics[width=10cm]{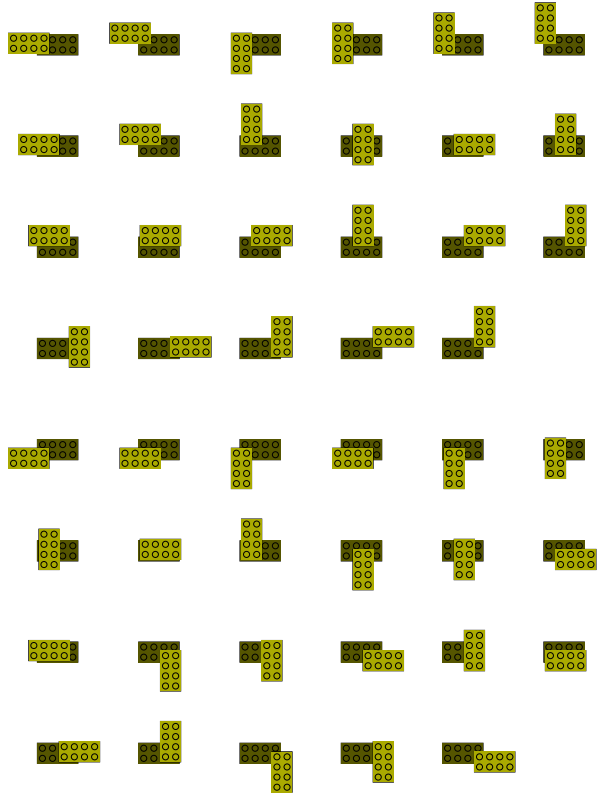}
\end{center}\caption{Even distribution}
\label{oneononeii}\end{figure}

On the other side of blocks $3,\dots,n-1$, we can do even better, as
outlined above. We leave to the reader to check that the  30 positions may be distributed evenly over 5
studs. Hence,
\eqref{sixbounds} is established, and the remaining claim follows by
Stirling's formula as in Theorem \ref{genest}.
\end{demo}

% \begin{figure}
%\begin{center}
%\includegraphics[width=5cm]{oneoneiibackc}\qquad
%\includegraphics[width=5cm]{oneoneiibackm}\qquad
%\end{center}\caption{Even distribution, previously used side}\label{oneononeiiback}
%\end{figure}

It turns out -- somewhat counterintuitively? -- that uneven
distributions of the positions give slightly better estimates than
what we obtained above. We have not carried out a systematic analysis
and can not claim that the distribution leading to Theorem
\ref{unevenest} is optimal, 
but trial and error with the following proposition make us believe
that there is only marginal room for improvement by this method.

If $(a_1,\dots,a_8)$ and  $(b_1,\dots,b_8)$ are tuples of integers, we
write $(a_1,\dots,a_8)\leq(b_1,\dots,b_8)$ if there is a permutation
$\sigma$ of $\{1,\dots,8\}$ with the property that $a_{\sigma(i)}\leq b_i$
for each $i$.

The methods leading to the following result are surely known.

\begin{propo}{alldist}
Let $\mathcal{P}_i$ and $\mathcal{P}^j_i$ be partitions of the sets of positions as
outlined above, and assume that
\begin{gather*}
(\#\mathcal{P}_1,\dots,\#\mathcal{P}_8)\leq(a_1,\dots,a_8)\\
(\#\mathcal{P}^j_1,\dots,\#\mathcal{P}^j_8)\leq(b_1,\dots,b_8),\qquad j\in\{1,\dots,8\}.
\end{gather*}
With
\[
{P}_0(y)=(y+a_1)\cdots(y+a_8)\qquad  P(y)={P}_0(y)(y+b_1)\cdots(y+b_8)
\]
we have that $\aaa_n$ is dominated by the coefficient of $y^{15n-31}$
in $({P}_0(y))^2(P(y))^{n-3}$ and that 
\[
h_\tbfm\leq \frac{P(x_0)}{x_0^{15}}
\]
where $x_0$ is the largest real root of
\[
Q(x)=15P(x)-xP'(x)
\]
\end{propo}

%GIVE PROOF OR NOT?

With the even distribution described above we get $x_0=72$, which,
since $P(72)/72^{15}=6\cdot13^{13}/12^{12}$  is consistent with Theorem
\ref{sixest}. Using that in fact
$(\#\mathcal{P}_1,\dots,\#\mathcal{P}_8)\leq(5,5,6,6,6,6,6,6)$ we may improve the estimate
on $h_\tbfm$ sligthly to $198.57$. 

However, we can do even better with very uneven
distributions:

\begin{theor}{unevenest}
$h_\tbfm\leq 191.35$
\end{theor}
\begin{demo}
 There exist partitions with
\begin{gather*}
(\#\mathcal{P}_1,\dots,\#\mathcal{P}_8)\leq(16,15,7,5,2,1,0,0)\\
(\#\mathcal{P}^j_1,\dots,\#\mathcal{P}^j_8)\leq(15,7,4,3,1,0,0,0),\qquad j\in\{1,\dots,8\}.
\end{gather*}
so by Proposition \ref{alldist} we are lead to consider
\[
P(x)=-x^5(x+15)(x+7)(x+1)R(x)
\]
where $R(x)$ is the polynomial
\[
x^8-23x^7-2056x^6-38700x^5-332657x^4-1504645x^3-3645736x^2-4392600x-2016000
\]
which has a largest real root which is approximately $65.05$. The estimate
follows by Proposition \ref{alldist}.
\end{demo}

\section{Improved lower bounds}

It follows from Theorem \ref{genest} that $h_\tbfm\geq 46$. We shall
in this section improve this estimate to  $h_\tbfm> 78.32$.

We let $\BBB_n\subset\AAA_{n+1}$ denote the set of \textsf{LEGO} configurations as above
consisting of $n+1$ blocks and such that both the top and the bottom
layer consists of a single block. Setting $\bbb_n=\#\BBB_n$ we then
clearly have
\[
\aaa_n\leq \bbb_n\leq \aaa_{n+1}
\]
and hence, as in Lemma \ref{sameentro},
\begin{eqnarray}\label{sameent}
h_\tbfm=\exp\left(\lim_{n\longrightarrow\infty}\frac{\log\bbb_n}{n}\right).
\end{eqnarray}

We say that a configuration $c$ in $\BBB_n$ has a
\emph{bottleneck} at height $z\in\NN$ if $c$ has exactly one
block in the layer $\RR^2\times [z;z+1]$. By convention the
top and bottom blocks are  not bottlenecks. This ensures
that removal of a bottleneck decomposes $c$ into two configurations
$c_0'$ and $c_0''$ one of which, say $c_0'$, contains the bottom block
of $c$. Re-inserting the removed block in $c_0'$ yields a
configuration $c'$ in some $\BBB_n$ with the inserted block as the top
block. Re-inserting the removed block into $c_0''$ yields, after a
translation, a configuration $c''$ in $\BBB_{n-m}$ with the inserted
block (translated) as the bottom block. Evidently, we can reconstruct
$c$ in a unique fashion from $(c',c'')$. Repeating this decomposition
procedure we conclude that any configuration in $\BBB_n$ with exactly
$k\geq 0$ bottlenecks can in a unique way be decomposed into a sequence
$(c^{(1)},\dots,c^{(k+1)})$ of configurations such that
$c^{(i)}\in\BBB_{m_i}$ has no bottlenecks and $m_1+\dots+m_{k+1}=n$.
 \begin{figure}
\begin{center}
\includegraphics[width=7cm,angle=90]{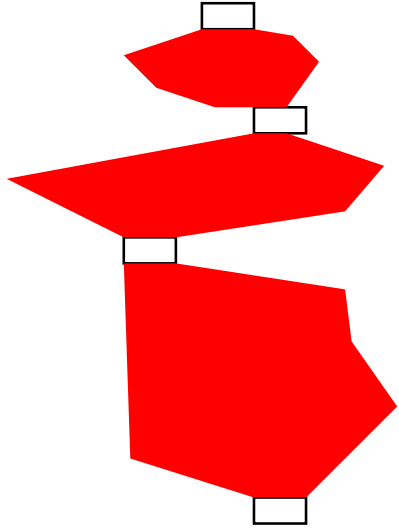}
\end{center}\caption{Two bottlenecks}\label{bneck}
\end{figure}

Letting $\CCC_n$ denote the subset of $\BBB_n$ consisting of
configurations without bottlenecks we obtain in this way a one-to-one
correspondence between elements of $\BBB_n$ and those of
\begin{equation}\label{bdecompintoc}
\bigcup_{k=
  0}^\infty\left[\bigcup_{m_1+\dots+m_{k+1}=n}\CCC_{m_1}\times\cdots\times
\CCC_{m_{k+1}}\right]
\end{equation}
Let  now
\[
\ccc_n=\#\CCC_n
\]
and let $\psi$ and $\psi_0$ denote the generating functions
\begin{equation}\label{genfuncts}
\psi(z)=\sum_{n=1}^\infty\bbb_nz^n\qquad
\psi_0(z)=\sum_{n=1}^\infty\ccc_nz^n
\end{equation}
It follows from \eqref{bdecompintoc} that
\begin{equation}\label{psipisn}
\psi(z)=\sum_{i=1}^\infty(\psi_0(z))^i=\frac{\psi_0(z)}{1-\psi_0(z)}.
\end{equation}
From the definition of $\psi(z)$ and $h_\tbfm$ it follows that $\psi$
is analytic in the disc
\[
D=\{z\mid|z|<(h_\tbfm)^{-1}\}
\]
and, since $b_n>0$, that $\psi$ is non-analytic at
$z=(h_\tbfm)^{-1}$. From \eqref{psipisn} we hence conclude that
\begin{equation}\label{lessone}
|\psi_0(z)|<1\qquad\text{ for }|z|<(h_\tbfm)^{-1}
\end{equation}
In particular, we get
\[
\ccc_1(h_\tbfm)^{-1}+\cdots+\ccc_n(h_\tbfm)^{-n}\leq 1
\]
which gives our claimed lower bounds on $h_\tbfm$, depending on the
number of terms $n$ on the lefthand side. 

We shall describe in detail how to get the first order of improvement
of the estimate from Theorem \ref{genest}. As evidently $\ccc_2=0$ we
turn to $\ccc_3$ for this.

The configurations contributing to $\ccc_3$ have one bottom block, one
top block. and two blocks in between. The number of ways of placing
two blocks on top of the bottom block is rather easily seen to be
$480$, so the number of configurations where the two middle blocks are
both attached to the bottom block is $2\cdot46\cdot 480-4730$, where
$4730$ is the number of configurations where the middle blocks are
both attached to the top block as well as to the bottom block. The
remaining configurations are those where the middle blocks are both
attached to the top block but only one of them to the bottom. This
number is seen to be $2\cdot46\cdot480-2\cdot4730$. Thus we have
\[
\ccc_3=4\cdot 46\cdot 480-3\cdot 4730=74130
\]
and $P_3(h_\tbfm^{-1})\leq 1$ where
\[
P_3(x)=46x+{74130}x^3
\]
This gives
\[
h_\tbfm\geq 64.06
\]
which can be improved as follows.

\begin{theor}{bestlo} $h_\tbfm> 78.32$
\end{theor}
\begin{demo}
Computer-aided computations give
\begin{eqnarray*}
\ccc_4&=&867346\\
\ccc_5&=&318434429\\
\ccc_6&=&18335373238
\end{eqnarray*}
so we have that $P_6(h_\tbfm^{-1})\leq 1$ where
\[
P_6(x)=46x+{74130}x^3
+{867346}x^4
+{318434429}x^5
+{18335373238}x^6
\]
which gives $h_\tbfm> 76.67$.

To improve the estimate we prove that
\begin{equation}\label{cest}
\ccc_{n+2}\geq 1248\ccc_{n}
\end{equation}
for $n\geq 6$. To see this, we devise 1248 different ways to construct an
element $c'$ of $\CCC_{n+2}$ from an element $c$ of $\CCC_n$, in such a way
that the original configuration can be recovered from the resulting
one.

Of these 1248 configurations, 480 are gotten from a fixed
configuration $d$ of two \bls sitting on one base \bl $b$ by identifying the base
block of $c$ with the block at the second level of $d$ which meets the
stud of lowest index on $b$ according to the enumeration of Figure
\ref{numerate}. We get the configuration $c'$ by rotating 90$^\circ$,
if necessary. The remaining 768 configurations are gotten by placing
one \bl underneath the base block of $c$, and placing one more block
at the level of this original base block, such that these two added
blocks do not meet. A computer search shows that there is always at
least this number of ways to do so since there are at least two blocks
at level $1$ of $c$ but only one at level $0$. We get the configuration $c'$ by
translating the configuration upwards and rotating 90$^\circ$,
if necessary.

To reconstruct $c$ from $c'$, one first sees how many blocks are
attached to the base block of $c'$. If there are two, $c$ is gotten by
discarding the base block of $c'$ and the block  at the next level
sitting at the highest index of it, translating down and rotating $270^\circ$, if
necessary. If there is only one, $c$ is gotten by discarding the base
block of $c'$ and the block at the next level which does not meet that
block, translating down and rotating $270^\circ$, if
necessary.

By repeated application of \eqref{cest} we get $\ccc_{6+2k}\geq
1248^k\ccc_6$  and $\ccc_{5+2k}\geq
1248^k\ccc_5$, so that  $r(h_{\tbfm})\leq 1$ with
\[
r(x)=P_6(x)+\frac{\ccc_5x^7+\ccc_6x^8}{1/1248-x^2},
\]
leading to the stated lower bound.
\end{demo}

\section{Concluding remarks}
We do not at present have the software nor the computer power to
perform numerical experiments to get a good idea of the true value of
$h_{\tbfm}$. Our best guess, based mainly on data achieved by Abrahamsen on the presumably closely
related case of $1\times 2$-blocks would be that the number is rather close
to $100$.

\begin{center}\sc Department of Mathematics\\University
  of Copenhagen\\Universitetsparken 5\\DK-2100 Copenhagen Ø\\Denmark
\end{center}
\begin{thebibliography}{1}
\bibitem{ma:data} Mikkel Abrahamsen, \emph{\textsf{LEGO} counting results}. Odsherreds
  Gymnasium. 
\bibitem{jkc} Jørgen Kirk Christiansen, \emph{"Taljonglering med klodser -- eller talrige klodser}. Klodshans (\textsf{LEGO} Company newsletter), 1974.
\bibitem{ulb} Kjeld Kirk Christiansen, \emph{The ultimate \textsf{LEGO}
  book}. Dorling Kindersley, 1999.
\bibitem{lcp} \textsf{LEGO} Company Profile 2004, \lego, \url{http://www.lego.com/info/pdf/compprofileeng.pdf}.
\bibitem{llife} Trine Nissen, \emph{From 102 to 915 million
    combinations} \textsf{LEGOLife} 5, Company newsletter, 2004. Available from
\url{http://www.math.ku.dk/~eilers/lego.html}.
\end{thebibliography}
\end{document}